\documentstyle{article}
\begin{document}
\begin{center}
{{\LARGE Stability of Polytopic Polynomial Matrices\footnote[1]
 { Supported by National Key Project and National Natural Science Foundation of China(69925307). Email: longwang@mech.pku.edu.cn}
}}
\end{center}

\vskip 0.6cm
\centerline{Long Wang \hspace{2cm} Zhizhen Wang \hspace{2cm} Wensheng Yu}
\vskip 6pt
\centerline{\small{Center for Systems and Control, Department of Mechanics and Engineering Science }}
\centerline{\small{Peking University, Beijing 100871, CHINA }}

\vskip 0.6cm

\vskip 6pt
\begin{minipage}[t]{14cm}{
\noindent{Abstract:}
This paper gives a necessary and sufficient condition for robust $D$-stability of 
Polytopic Polynomial Matrices.
}
\end{minipage}

\section{Introduction}
\par \indent

Motivated by the seminal theorem of Kharitonov on robust stability of interval 
polynomials\cite{Khar, Khar1}, a number of papers on robustness analysis of uncertain
systems have been published in the past few years\cite{Holl, Bart, Fu, wang1, 
wang2, Barmish, Chap, wang3}. Kharitonov's theorem states that the Hurwitz stability of
a real (or complex) interval polynomial family can be guaranteed by the Hurwitz 
stability of four (or eight) prescribed critical vertex polynomials in this family.
This result is significant since it reduces checking stability of infinitely many 
polynomials to checking stability of finitely many polynomials, and the number of
critical vertex polynomials need to be checked is independent of the order of the 
polynomial family. An important extension of Kharitonov's theorem is the edge theorem
discovered by Bartlett, Hollot and Huang\cite{Bart}. The edge theorem states that
the stability of a polytope of polynomials can be guaranteed by the stability of its
one-dimensional exposed edge polynomials. The significance of the edge theorem is 
that it allows some (affine) dependency among polynomial coefficients, and applies to
more general stability regions, e.g., unit circle, left sector, shifted half plane,
hyperbola region, etc. When the dependency among polynomial coefficients is nonlinear,
however, Ackermann shows that checking a subset of a polynomial family generally can 
not guarantee the stability of the entire family\cite{Ack1, Ack2, Ack3}.

Parallel to this line of research, robust stability of uncertain matrices has also received considerable attention.
Bialas 'proved' that for robust Hurwitz stability of an interval matrix, it suffices to check all vertices\cite{Bialas}. But Barmish and Hollot
gave a counter-example to show that Bialas's claim is incorrect\cite{Barmish1}. Kokame and Mori considered
Hurwitz stability of an interval polynomial matrix, and by using some result in signal processing theory, established
a necessary and sufficient condition for robust stability\cite{Mori}.

This paper studies robust $D$-stability of polytopic polynomial matrices, i.e., matrices with entries being
polytopes of polynomials. We give a necessary and sufficient condition for robust $D$-stability of 
Polytopic Polynomial Matrices, namely, the stability of a subset of this family guarantees the stability
of the entire family.

\vskip -1in
\section{\bf{ Preliminaries}}
\par \indent

{
{\bf Definition 1}\quad An interval polynomial matrix
$A=(p_{ij})_{n\times n}$ is a matrix whose entries $p_{ij}$ are
interval polynomials, i.e., $p_{ij}=q_{ij}^0+q_{ij}^1 s+\dots+q_{ij}^m
s^m, \quad q_{ij}^k \in[\underline q_{ij}^k, \overline q_{ij}^k],
\quad k=0, \dots, m$, where $k$ stands for superscript.
A polytopic polynomial matrix $A=(p_{ij})_{n\times n}$ is a matrix
whose entries $p_{ij}$ are polytopic polynomials, i.e., $p_{ij}=
\sum_{k=1}^m \lambda_{ij}^k p_{ij}^k, \quad
\lambda_{ij}^k \geq 0, \quad \sum_{k=1}^m \lambda_{ij}^k=1, \quad
i, j=1, \dots, n$, where $p_{ij}^k$ are fixed polynomials.

{\bf Definition 2}\quad  Suppose $\Omega \subset{\bf R}^{n+1}$ is an
$m$-dimensional polytope. Its supporting plane $H$ is defined as an
$n$-dimensional affine set, satisfying $\Omega \cap H \not=\emptyset$,
and all points of $\Omega$ lie on the same side of $H$; Its exposed
set is defined as the intersection of $\Omega$ and its
supporting plane $H$; Its exposed edge set is defined as the
one-dimensional exposed set.

An $n$-th order polynomial can be regarded as a point in its
$(n+1)$-dimensional coefficient space.

{\bf Definition 3}\quad Given an open region $D$ in the complex plane,
the polynomial matrix $A$ is said to be $D$-stable, if all roots of
$detA=0$ lie within $D$; A polynomial matrix set ${\cal A}$ is said
to be $D$-stable, if every member in ${\cal A}$ is $D$-stable. When
$D$ is taken as the open left half of the complex plane,
$D$-stable is also called Hurwitz stable.

In fact, Edge Theorem holds for more general stability regions. For
simplicity, we only consider simply-connected stability regions
in this paper.

{\bf Definition 4}\quad Given an interval polynomial set ${\cal F}(s)=
\{\sum_{i=0}^m q_i s^i, \quad q_i \in [\underline q_i, \overline q_i]
\}$, its Kharitonov vertex set is
 $K_{\cal F}^0=\{f_k^1, f_k^2, f_k^3, f_k^4\}$, and $E_{\cal F}^0=
\{\lambda f_k^s+(1-\lambda) f_k^t, \quad
(s,t)\in \{(1,2), (2,4), (4,3), (3,1)\}, \lambda \in [0,1]\}$ is
called its Kharitonov exposed edge set,  where
  $$
  {
  \begin{array}{ll}
  f_k^1=\underline q_0+\underline q_1 s+\overline q_2 s^2+
   \dots   &
  f_k^2=\underline q_0+\overline q_1 s+\overline q_2 s^2+
   \dots \\
  f_k^3=\overline q_0+\underline q_1 s+\underline q_2 s^2+
   \dots   &
  f_k^4=\overline q_0+\overline q_1 s+\underline q_2 s^2+
   \dots
  \end{array}
  }
  $$
Consider the polytopic polynomial sets
\begin{equation}
\begin{array}{c}
{\cal P}_{ij}=\left \{\sum_{k=1}^m \lambda_{ij}^k p_{ij}^k:
\lambda_{ij}^k \geq 0 , \sum_{k=1}^m \lambda_{ij}^k=1 \right \}
(i, j=1, \dots, n)\\  p_{ij}^k \mbox{ are fixed polynomials, }
k=1, \dots, m.
\end{array}
\end{equation}
Their vertex sets are
$$
{ K_{ij}=\{p_{ij}^k \quad k=1, \dots, m \} \quad   (i, j=1, \dots, n)}
$$
and by definition, their exposed edge sets are contained in
$$
{ E_{ij}=\{ \lambda p_{ij}^s+(1-\lambda)p_{ij}^t,\; s,t=1, \dots, m \}
 \quad (i, j=1, \dots, n)}
$$
Let
${\cal A} = \{(p_{ij})_{n \times n}: p_{ij} \in {\cal P}_{ij},\;
i, j=1, \dots, n \}$, and let $P_n^n$ be the set
of all permutations of $1, 2, \dots, n$.

{\bf Definition 5}\quad  Define $\epsilon_{\cal A}$ as
\begin{equation}
\left \{(p_{ij})_{n \times n}:
\begin{array}{l}
p_{s l_s}\in E_{s l_s} , (l_1, \dots, l_n )\in P_n^n , s=1, \dots, n \\
p_{s i_s} \in K_{s i_s} , i_s=1, \dots, l_s -1, l_s +1, \dots, n
\end{array}
\right \}
\end{equation}
It is easy to see that, $\epsilon_{\cal A}$ is produced by taking only
one entry from its exposed edge set in every row/column
and all other entries from their vertex sets in ${\cal A}$.}
\vskip -2in
 \section{\bf  Main Results}
 \subsection{\bf Polytopic Polynomial Matrices}
 \par \indent

Consider the polytopic polynomial matrix set
\begin{equation}
{\displaystyle {\cal A}=\{
(p_{ij})_{n \times n} : \mbox{ where }
 p_{ij} \in {\cal P}_{ij},\; i,j=1, \dots, n \}}
\end{equation}
 Suppose $\forall A \in {\cal A}, deg(detA)=const.$

{\bf Theorem 1}\quad ${\cal A}$ is $D$-stable if and only if
$\epsilon_{\cal A}$ is $D$-stable.

Proof: Necessity is obvious. To prove sufficiency, suppose
$\epsilon_{\cal A}$ is $D$-stable. Let
$$
 A=\left(
\begin{array}{rcl}
\renewcommand\arraystretch{0.2}
p_{11}& \dots & p_{1n}\\
\dots & \dots & \dots \\
p_{n1}& \dots & p_{nn}
\end{array}
\right)
$$
By using Laplace formula on the first column, we have
$$
detA=p_{11}M_{11}+\dots+p_{n1}M_{n1}
$$
Let
$$
{\cal T}=\left \{
\left(
\begin{array}{rccl}
\arraycolsep=.01mm
\renewcommand\arraystretch{0.2}
p_{11}^{*} & p_{12} & \dots & p_{1n} \\
\dots & \dots & \dots & \dots \\
p_{n1}^{*} & p_{n2}& \dots & p_{nn}
\end{array}
\right) ,
 \begin{array}{l}
\mbox{and } p_{i1}^{*} \in {\cal P}_{i1} ;  i=1, \dots, n \\
p_{ij} \mbox{ are entries of } A \\
j=2, \dots, n
\end{array}
\right \}
$$
It is easy to see that $A \in {\cal T}$ and $\forall \; T \in {\cal T}$
$$
detT=p_{11}^{*}M_{11}+\dots+p_{n1}^{*}M_{n1}
$$
Apparently, $detT$ is an affine function of $p_{11}^{*}, \dots,
p_{n1}^{*}$. By Edge Theorem
$$
detT\mbox{ is }D \mbox{-stable}\Leftrightarrow \mbox{the edge set of }
detT\mbox{ is }D \mbox{-stable.}
$$
The edge set of $detT$ is
$$ \left \{
E_{i1}M_{i1}+\sum_{i\not=j=1}^{n}K_{j1}M_{j1},\; i=1, \dots, n
\right \}.
$$
The corresponding matrix collection is
$${\cal A}_1=\left \{\left(
\begin{array}{c}
\renewcommand\arraystretch{0.2}
q_{11}\\
\dots \\
q_{i1}\\
\dots \\
q_{n1}
\end{array}
(p_{st})_{n \times (n-1)}
\right),
\begin{array}{l}
\renewcommand\arraystretch{0.2}
q_{i1} \in E_{i1} ; i \in \{1, \dots, n \}\\
q_{k1}\in K_{k1} ; t=2, \dots, n\\
p_{kt} ,\; p_{it}\mbox{ are entries of } A \\
k=1, \dots,i-1,i+1,\dots, n
\end{array}
\right \}$$
In this case, ${\cal T}$ is $D$-stable $\Leftrightarrow {\cal A}_1$ is
$D$-stable. Moreover, $\forall A_1\in {\cal A}_1$, there exists
$i \in \{1, \dots, n \}, q_{i1} \in
E_{i1}, q_{k1} \in K_{k1}, k=1, \dots, i-1, i+1, \dots, n $ such that
$$
A_1=\left(
\begin{array}{cc}
\renewcommand\arraystretch{0.2}
q_{11}&p_{12} \\
\dots &\dots \\
q_{i1}&p_{i2} \\
\dots &\dots \\
q_{n1}&p_{n2}
\end{array}
(p_{st})_{n \times (n-2)}
\right)
$$
Again, by using Laplace formula on the second column, we have
$$
detA_1=p_{12}M_{12}+\dots+p_{n2}M_{n2}
$$
Set
$$
{\cal B}=\left \{
\left(
\begin{array}{ccccc}
\arraycolsep=.01mm
\renewcommand\arraystretch{0.2}
q_{11}&p_{12}^{*}&p_{13}&\dots&p_{1n}\\
\dots&\dots&\dots&\dots&\dots\\
q_{n1}&p_{n2}^{*}&p_{n3}&\dots&p_{nn}
\end{array}
\right)
\begin{array}{l}
\arraycolsep=.01mm
\renewcommand\arraystretch{0.2}
q_{ij} \mbox{ are entries of } A_1\\
p_{ij} \mbox{ are entries of } A \\
p_{ij}^{*}  \in {\cal P}_{ij}\\
i,j=1, \dots, n
\end{array}
\right\} $$ Let ${\cal T}_1=\bigcup_{A_1 \in {\cal A}_1} {\cal
B}$, then, its edge set ${\cal A}_2$ is $$ \left \{ \left(
\begin{array}{ccl}
\arraycolsep=.01mm
\renewcommand\arraystretch{0.2}
q_{11}& q_{12}& \\
\dots & \dots & \\
q_{i1}& \dots & \\
\dots & \dots &
\begin{array}[b]{l}
(p_{st}) \\
\end{array}\\
\dots & q_{j2}& \\
\dots & \dots & \\
q_{n1}& q_{n2}&
\end{array}
\right)
\begin{array}{l}
j \in  \{1, \dots, n \}; l=1, \dots, n \\
q_{lm} \in \left\{
\begin{array}{l}
 K_{lm};
 \begin{array}{l}
 m=1, l \not=i\\
 m=2, l \not=j
 \end{array}\\
 E_{lm};
 \begin{array}{l}
 m=1, l=i\\
 m=2, l=j
 \end{array}
 \end{array}
  \right.\\
p_{lt} \mbox{ are entries of } A \\
t=3, \dots, n
\end{array}
 \right\}
 $$
 and ${\cal A}_1 \subset {\cal T}_1$. By definition and Edge Theorem
$$
{\cal T}_1\mbox{ is }D \mbox{-stable }\Leftrightarrow {\cal A}_2
\mbox{ is }D \mbox{-stable.}
$$
By repeating the process above, we have
$$
\begin{array}{l}
{\cal A}_2, \dots,  {\cal A}_{n-1}, {\cal A}_n;\quad
{\cal T}_2, \dots, {\cal T}_{n-1} \\
{\cal T}_k \mbox{ is }D \mbox{-stable}\Leftrightarrow {\cal A}_{k+1} \mbox{ is }
D \mbox{-stable}\\
{\cal A}_k \mbox{ is }D \mbox{-stable} \Leftarrow {\cal T}_k
\mbox{ is } D \mbox{-stable }
\quad k=2, \dots, n-1
\end{array}
$$
where ${\cal A}_k$ is the collection of
$$\left \{\left(
\begin{array}{rccc}
\arraycolsep=.01mm
\renewcommand\arraystretch{0.2}
q_{11}& \dots & q_{1k}& \\
\dots & \dots & \dots & \\
q_{i_1 1}& \dots & q_{i_1 k} & \\
\dots & \dots & \dots &(p_{st})\\
q_{i_k 1}& \dots & q_{i_k k}& \\
\dots & \dots & \dots & \\
q_{n1}& \dots & q_{nk}&
\end{array}
\right)_{n \times n}
\begin{array}{l}
q_{lt} \in \left\{
\begin{array}{l}
 E_{i_t t} , l=i_t \\
 K_{i_t t} , l \not=i_t
\end{array}
\right.\\
i_1, \dots, i_k \in  \{1, \dots, n \} \\
t=1, \dots, k \\
l=1, \dots, n \\
p_{ls} \mbox{ are entries of } A \\
s=k+1, \dots, n
\end{array}
\right \} $$ Thus, for each element of $ {\cal A}_k$, its entries
have the following characteristics: for the first $k$ columns, all
entries of each column belong to their vertex sets except that one
entry belongs to its exposed edge set, and the entries of the
remaining $n-k$ columns are the corresponding entries of $ A$.
Hence, $\forall A_n \in {\cal A}_n$, we have $$ {\cal
A}_n=\left\{\left(
\begin{array}{rcl}
\arraycolsep=.01mm
\renewcommand\arraystretch{0.2}
q_{11}& \dots & q_{1n}\\
\dots & \dots & \dots \\
q_{i_1 1}& \dots & q_{i_1 n}\\
\dots & \dots & \dots \\
q_{i_n 1}& \dots & q_{i_n n}\\
\dots & \dots & \dots \\
q_{n1}& \dots & q_{nn}
\end{array}
\right)
\begin{array}{l}
q_{lt} \in \left\{
\begin{array}{l}
 E_{i_t t} , l=i_t \\
 K_{i_t t} , l \not=i_t
\end{array}
\right.\\
i_1, \dots, i_n \in  \{1, \dots, n \} \\
t=1, \dots, n \\
l=1, \dots, n
\end{array}
\right \} $$
If $i_s=i_t$ for some pair $i_s, i_t$, without loss
of generality, suppose $i_1=i_2=1$, namely $$A_n=\left(
\begin{array}{rccl}
\arraycolsep=.01mm
\renewcommand\arraystretch{0.2}
q_{11}& q_{12}& \dots & q_{1n}\\
\dots & \dots & \dots & \dots \\
q_{n1}& q_{n2}& \dots & q_{nn}
\end{array}
\right)
,\; q_{11} \in E_{11},\; q_{12} \in E_{12}
$$
By using Laplace formula on the first row of $A_n$, we have
$$
detA_n=q_{11}M_{11}+q_{12}M_{12}+\sum_{i=3}^n q_{1i}M_{1i}
$$
By Edge Theorem
$$
\begin{array}{ll}
\arraycolsep=.01mm
A_n \mbox{ is }D \mbox{-stable}& \Leftrightarrow
q_{11}M_{11}+q_{12}^{0}M_{12}+\sum_{i=3}^n q_{1i}M_{1i}  \mbox{ and }\\
 \end{array}
$$
{\parskip=-0.1in
$$
q_{11}^{0}M_{11}+q_{12}M_{12}+\sum_{i=3}^n q_{1i}M_{1i}\mbox{ are }D \mbox{-stable.}
$$}
The corresponding matrices are
$$
 \left(
\begin{array}{rccl}
 \renewcommand\arraystretch{0.2}
q_{11}^{0}& q_{12}& \dots & q_{1n}\\
\dots & \dots & \dots & \dots \\
q_{n1}& q_{n2}& \dots & q_{nn}
\end{array}
\right)
, q_{11}^{0} \in K_{11},
$$
$$
\left(
\begin{array}{rccl}
 \renewcommand\arraystretch{0.5}
q_{11}& q_{12}^{0}& \dots & q_{1n}\\
\dots & \dots & \dots & \dots \\
q_{n1}& q_{n2}& \dots & q_{nn}
\end{array}
\right)
,   q_{12}^{0} \in K_{12},
 $$
which belong to
$$ \left \{\left(
\begin{array}{rcl}
\arraycolsep=.04mm
\renewcommand\arraystretch{0.5}
q_{11}& \dots & q_{1n}\\
\dots & \dots & \dots \\
q_{i_1 1}& \dots & q_{i_1 n}\\
\dots & \dots & \dots \\
q_{i_n 1}& \dots & q_{i_n n} \\
\dots & \dots & \dots \\
q_{n1}& \dots & q_{nn}
\end{array}
\right),
\begin{array}{l}
q_{i_s s}\in E_{i_s s} , (i_1, \dots, i_n )\in P_n^n  \\
q_{ls} \in K_{ls} ,l \not=i_s \\
l=1, \dots, n \\
s=1, \dots, n
\end{array}
\right \}.
$$
So, $ {\cal A}_n \mbox{ is }D \mbox{-stable} \Leftrightarrow
\epsilon_{\cal A}\mbox{ is }D \mbox{-stable}$. Thus
$$
\begin{array}{l}
\renewcommand\arraystretch{0.5}
\epsilon_{\cal A}\mbox{ is }D \mbox{-stable} \Leftrightarrow
{\cal A}_n \mbox{ is }D\-- \mbox{stable}
 \Leftrightarrow {\cal T}_{n-1} \mbox{ is } D \mbox{-stable}\\
 \Rightarrow {\cal A}_{n-1} \mbox{ is } D \mbox{-stable}
 \dots
 \Leftrightarrow {\cal T}_i \mbox{ is } D \mbox{-stable}\\
 \Rightarrow {\cal A}_i \mbox{ is } D \mbox{-stable}
 \Leftrightarrow {\cal T} \mbox{ is } D \mbox{-stable}
 \Rightarrow A \mbox{ is } D \mbox{-stable.}
\end{array}
$$
That is to say, $\epsilon_{\cal A}\mbox{ is }D \mbox{-stable}
\Rightarrow \forall A \in {\cal A}, A \mbox{ is }D \mbox{-stable}$,
namely, ${\cal A}$ is $D$-stable. This completes the proof.

{\bf Remark 1}\quad When $m=2$, i.e. $A=(p_{ij})_{n \times n}$,
where
\begin{equation}
\begin{array}{l}
p_{ij}=p_{ij}^0+\lambda_{ij} p_{ij}^1 \quad \lambda_{ij} \in [0, 1]\\
p_{ij}^0, p_{ij}^1 \mbox{ are fixed polynomials } \quad i, j=1, \dots, n
\end{array}
\end{equation}
Let
\begin{equation}
{\cal B}_1=\{(p_{ij})_{n \times n} \quad p_{ij}\mbox{ satisfies (4)}\}
\end{equation}
In this case, the vertex sets of $p_{ij}$ are $\{p_{ij}^0,p_{ij}^0+
p_{ij}^1 \}$, and their exposed sets are exactly themselves, namely,
$\{p_{ij}^0+\lambda_{ij}p_{ij}^1 \} $.
The corresponding conclusion has more concise form, this is due to the
simplification of the edge set of ${\cal B}_1$
$$
\epsilon_{{\cal B}_1}=\left\{
(p_{ij})_{n \times n}
\begin{array}{l}
(l_1, \dots, l_n) \in P_n^n ; s=1, \dots, n \\
\lambda_{si} \in \left \{
\begin{array}{l}
 [0, 1] \quad i=l_s\\
 \{0, 1\} \quad i\not=l_s
\end{array}
\right.
\end{array}
\right \}
$$
\subsection{\bf  Interval Polynomial Matrices}
 \par \indent

  Consider the subset of $ R^{n \times n}(s)$
\begin{equation}
 {\cal B}_2=\{(p_{ij})_{n \times n}, \; p_{ij}\mbox{ are interval
polynomials }\}
\end{equation}
where $R^{n \times n}(s)$ is the collection of $n \times n $ polynomial
matrices. Assume $\epsilon_{{\cal B}_2}$ is $Hurwitz$ stable.

Consider $detA$, using $Laplace$ formula on first column, we have
$$
detA=p_{11}M_{11}+\dots+p_{n1}M_{n1}\; .
$$
Similar to the proof of theorem 1, and by resort to  the Generalized Kharitonov
Theorem, we have

 {\bf Theorem 2}\quad If $\forall A \in{\cal B}_2, deg(detA)=m$. Then
${\cal B}_2$ is $Hurwitz$ stable $\Leftrightarrow \epsilon_{{\cal B}_2}$
is $Hurwitz$ stable, where $\epsilon_{{\cal B}_2}$ is
$$
\left \{(p_{ij})_{n \times n}:
\begin{array}{l}
p_{il_i}\in E_{i l_i}^0 ;  i=1, \dots, n \\
p_{si_s} \in K_{si_s}^0 ;  l_s \not=i_s=1, \dots, n ;\\
(l_1, \dots, l_n)\in P^n_n ;s=1, \dots, n
\end{array}
\right \}
$$
and $E_{i l_i}^0, K_{si_s}^0 $ are defined in Definition 4.

{\bf Remark 2}\quad Theorems 1 and 2 can be viewed as a generalization of 
the Edge Theorem and Kharitonov Theorem to MIMO case. Stability test of 
the entire family is reduced to a critical low-dimensional subset. No 
extra lemma from signal processing is needed in our proof. Furthermore, 
our results can be easily extended to polynomial matrices with complex 
coefficients.


\begin{thebibliography}{19}
 


\bibitem{Khar}V.L.Kharitonov. Asymptotic stability of an equilibrium
  position of a family of systems of linear differential
  equations, Differential'nye Uravneniya, vol.14, 2086-2088, 1978.

\bibitem{Khar1}V.L.Kharitonov. The Routh-Hurwitz problem for families of polynomials and quasipolynomials, Izvetiy Akademii Nauk Kazakhskoi SSR, Seria fizikomatematicheskaia, vol.26, 69-79, 1979.

\bibitem{Holl}C.V.Hollot and R.Tempo. On the Nyquist envelope of an interval plant family, IEEE Trans. on Automatic Control, vol.39, 391-396, 1994.
   
\bibitem{Bart}A.C.Bartlett, C.V.Hollot and L.Huang. Root locations of an entire polytope of polynomials: It suffices to check the edges, Mathematics of Control, Signals, and Systems, vol.1, 61-71, 1988.
  
\bibitem{Fu}M.Fu and B.R.Barmish. Polytope of polynomials with zeros in a prescribed set, IEEE Trans. on Automatic Control, vol.34, 544-546, 1989.

\bibitem{wang1}L. Wang and L. Huang. Vertex results for uncertain systems, Int. J. Systems Science, vol.25, 541-549, 1994.

\bibitem{wang2}L. Wang and L. Huang. Extreme point results for strict positive realness of transfer function families, Systems Science and Mathematical Sciences, vol.7, 371-378, 1994.
 
\bibitem{Barmish}B. R. Barmish, C. V. Hollot, F. J. Kraus and R. Tempo. Extreme point results for robust stabilization of interval plants with first order compensators, IEEE Trans. on Automatic Control, vol.37, 707-714, 1992. 

\bibitem{Chap}H. Chapellat, M. Dahleh and S. P. Bhattacharyya. On robust nonlinear stability of interval control systems, IEEE Trans. on Automatic Control, vol.36, 59-67, 1991.

\bibitem{wang3}L. Wang and L. Huang. Finite verification of strict positive realness of interval rational functions, Chinese Science Bulletin, vol.36, 262-264, 1991.
 
\bibitem{Ack1}J.Ackermann. Uncertainty structures and robust stability analysis, Proc. of European Control Conference, 2318-2327,1991.

\bibitem{Ack2}J.Ackermann. Does it suffice to check a subset of multilinear parameters in robustness analysis? IEEE Trans. on Automatic Control, vol.37, 487-488, 1992.

\bibitem{Ack3}J.Ackermann et al., Robust Control: Systems with Uncertain Physical Parameters, Springer-Verlag, Berlin, 1994.






\bibitem{Bialas}S.Bialas. A necessary and sufficient condition for the 
stability of interval matrices, Int. J. Control, vol.37, 717-722, 1983.

\bibitem{Barmish1}B.R.Barmish and C.V.Hollot. Counter-example to a recent result on the stability
of interval matrices by S. Bialas, Int. J. Control, vol.39, 1103-1104, 1984.

\bibitem{Mori}H.Kokame and T.Mori. A Kharitonov-like theorem for interval 
polynomial matrices, Systems and Control Letters, vol.16, 107-116, 1991.
 
\bibitem{Bhat}S.P. Bhattacharyya, H. Chapellat and L.H. Keel, Robust Control:
The Paramatric Approach, Prentice-hall, New Jersey, 1995.


















 
 
\bibitem{Ran}A.Rantzer. Stability conditions for polytopes of polynomials, IEEE Trans. on Automatic Control, vol.37, 79-89, 1992.



\bibitem{Zad}L.A.Zadeh and C.A.Desoer, Linear System Theory: A State Space Approach, McGraw-Hill, New York, 1963.

\bibitem{Barm}B.R.Barmish, New tools for robustness analysis, Proc. of IEEE Conf. on Decision and Control, 1-6, 1988.
 








\end{thebibliography}
\end{document}